# The Social and Work Structure of an Afterschool Math Club

Lawrence Smolinsky, Eugene Kennedy, Liuli Huang, and Andrew Alaniz
Louisiana State University

**Tables and figures are at the end of the article.**

## Abstract

The objective of this exploratory study was to examine the interpersonal dynamics of an afterschool math club for middle schoolers. Using social network analysis, two networks were identified and analyzed: (a) a network of friendship relationships and (b) a network of working relationships. The interconnections and correlations between friendship relationships, working relationships, and a student opinion survey were studied. We identified a core working group of students from within the network of working relations. This group acted as a central go-between for other members in the club. This core working group also expanded into the largest friendship group in the friendship network. A second group was formed from popular but aloof students who reported less impact from the club. Although there were working isolates, they were not found to be socially isolated. Students who were less popular tended to report a greater favorable impact from club participation than those who were more popular.

A growing body of research supports that participation in afterschool science, technology, engineering, and mathematics (STEM) programs and clubs has a positive impact on student engagement, achievement, and persistence in STEM majors and careers (Young et al., 2017; Krishnamurthi et al., 2014). These experiences afford students the flexibility to pursue their interests and engage in STEM disciplines in ways that are not possible in traditional classrooms where curriculum coverage and standardized tests often dictate the scope and pace of activities. In an afterschool STEM club, students may have (a) input into the activities of the club; (b) increased opportunities for student-to-student, student-to-teacher, and student-to-practicing STEM professional interactions; and (c) opportunities to receive mentoring and instruction tailored to their specific needs (Little et al., 2008). As a result, afterschool STEM club activities are often intrinsically rewarding and promote confidence, a sense of identity, self-efficacy, and enjoyment of STEM (Gmurczyk & Collins, 2010; Peterson, 2013). In the United States, STEM clubs are an increasingly important focus of public polity geared towards enhancing national competiveness in the STEM disciplines (National Research Council, 2009, 2015; STEM Education Coalition, 2016). According to Papanastasiou and Bottiger (2004), it is important to understand math clubs across international settings, reasoning that educational systems are intertwined with societal relationships and student opportunities outside of the normal school day. An educational system needs to be considered within its context (Robinsohn, 1992), including any informal education opportunities that may exist.

Despite the popularity of afterschool STEM clubs, little is known about the social and interpersonal dynamics of the children and youth who participate in the clubs. In general, students who voluntarily participate in these clubs tend to be somewhat homogeneous with respect to academic backgrounds and interests. However, research shows that some clubs are much more effective than others and that some students benefit more from participation in these experiences than their peers (e.g., Gottfried & Williams, 2013; Krishnamurthi et al., 2014; Young et al., 2017).

The National Research Council (2015) identified attributes of effective clubs. Similarly, several organizations and agencies have developed standards and guidelines for effective afterschool programs, including STEM clubs (e.g., Afterschool Alliance, 2011; Huang et al., 2008). Guidelines provided by these organizations tend to focus on structural features of these programs: administration, staff, interactions, activities, and so on. However, the social structure and network of interpersonal dynamics among participants in these clubs have rarely been examined.

The connections and interactions among members of groups can impact the experiences of participants, but the extent to which this typically occurs has not been systematically examined. For example, although it is commonly recognized that friendship networks affect the flow of information, the quality of collaborative work, and the co-construction of knowledge, little is known about how these processes unfold and operate in the informal context of an afterschool STEM club (Maroulis & Gomez, 2008). More specifically, within the context of mathematics learning, which is often characterized in terms of in-group and out-group contexts, club dynamics can be a significant factor in efforts to expand participation to include groups historically underrepresented in mathematics and other STEM disciplines. It is interesting that "being with friends" was the most selected reason to attend math club in a study within a Kansas City middle school (Papanastasiou & Bottiger, 2004). Patterns of friendships, the potential for isolates (individuals not connected to other members of the group), and access to talented peers with afterschool math clubs are issues that warrant attention.

The present study focused on a math club designed for low-income, primarily minority middle school students. We examined the social and working relationship among participating students as well as their attitudes toward mathematics, each other, and themselves. We employed a small sample size ($n$ = 19) and explored relationships among participants using network and graph theoretic techniques. Although these techniques have typically been used with large data sets, a small sample was used because these methods have not been previously applied to students in middle school math clubs. We anticipated that what was learned from the present study would allow us to refine the methods employed to use in the future with larger samples of students from multiple clubs.

**Research Questions**

In the current exploratory study, we used social network analysis to explore the following research questions:

- To what extent and in what ways are students' social relationships and working relationships dependent on each other?
- To what extent is a given student's position in the social structure of the club related to club experiences and/or characteristics of the club members?
- Are there distinct subgroups or communities within the club? If so, what are their characteristics?

## Literature Review

In this section we discuss the background of social learning in mathematics, math clubs, and social network analysis.

**Mathematics Learning as a Social Phenomenon**

The constructionist view of learning posits that individuals do not simply absorb knowledge, but rather knowledge is constructed through a dynamic process by which new ideas are tested and evaluated against existing ones, modified, and changed (e.g., Thompson, 2014). In the constructionist view, knowledge is personal, and the process of learning is inherently interactive and social. Mathematics education has, over the past several decades, increasingly integrated tenants of constructionism into teacher preparation programs and classroom instructional recommendations (Mueller et al., 2014).

Middle school students' collaborative work in mathematics classes has been examined within a variety of structured implementations (Nathan et al., 2007; Piccolo et al., 2008; Webb et al., 2006). Previous studies have focused on the importance of student-student interactions within their various structures. For example, Webb et al. (2006) pointed out that "…mathematics reform efforts consistently call for balanced and student-centered communication in which students take an active role in classroom discourse (Carpenter, Fennema, Peterson, Chiang, & Loef, 1989; Cobb & Baursfeld, 1995; Forman & Ansell, 2002)" (p. 64).

Discourses and *math talk* have allowed students to test and justify their ideas with their classmates and act as evaluators of the ideas and thoughts of others (e.g., Forman, 2003; Kosko, 2012). Investigations of classroom dialog and interactions among learners have increased significantly over the past years, with much of the focus being on how learners engage the process of co-construction of knowledge and how these processes are encouraged in some settings and not in others (Goos, 2004). Guidelines on productive discussions have typically focused on student behavior during discussion and the responses they offer to their peers (Razfar, 2012). Although these guidelines are important, largely absent from the literature have been guidelines relative to the desired social dynamics of the classroom.

**Afterschool Math Clubs**

Unlike traditional classrooms, afterschool math clubs are informal environments in which students often play a significant role in the selection of activities. Discussion, the co-construction of knowledge, and other collaborative experiences are common in afterschool clubs; these activities are often suggested as reasons these clubs are associated with positive impacts on student outcomes (Afterschool Alliance, 2011).

For example, Thompson (2009) found that math clubs had a positive impact on problem solving ability, and Gottfried and Williams (2013) reported positive effects on standardized test scores. The effects of math clubs may have extended beyond the K12 setting. Participation in afterschool clubs was positively related to retention in a STEM major and the likelihood of majoring in a math intensive major in post-secondary education (Sahin, 2013).

***Afterschool Math Clubs and Social Network Analysis***

Despite the growing body of research on the effects of afterschool math clubs, as with STEM clubs in general, there has been relatively little work on the social dynamics of clubs that may act to promote or inhibit their effectiveness. Social network analysis could offer insights into the inner workings of a club and inform efforts to promote club productivity. Social network analysis has been concerned with connections among people in groups or organizations. These connections could be based on friendships, interpersonal interactions, or other relational factors. Understanding these connections could provide insights into their impact on club outcomes. In education, social network analysis has had significant potential (Grunspan et al., 2014) and has been used to investigate learning communities in online learning environments (Cadima et al., 2012), collaborations among teachers (Lin et al., 2016), and peer effects (Santos et al., 2015).

**Social Network Analysis**

In mathematics, graph theory studies and classifies the complex structures that arise from pairwise relationships between individuals or objects called nodes (Thulasirama & Swamy, 1992). Simple mutual relationships between nodes give graphs (e.g., a and b are friends). Directed graphs may express relationships that are not symmetrical (e.g., a loves b). More complex relations may be expressed by directed weighted graphs, e.g., (1) a loves b, (2) a likes b, (3) a dislikes b, (4) a hates b. Labeled graphs are often call networks. Social network analysis investigates social relationships through the use of graph theory. For example, individuals may be designated as central, peripheral or isolated with respect to a social network. The applications and study of networks has expanded throughout the sciences and social science with the development of computing power (Newman, 2010).

## Method

We start with a brief overview with details described in the subsections. In 2013–2014, a similar outreach math club was conducted with another group of students (Kennedy & Smolinsky, 2016). That study focused on the impact of participation in the math club on students and the design features of the experience that were most effective at promoting engagement and positive reactions from students, but the recruitment and operation was like the present study. However, in the present study, students met in a classroom with tables rather than desks and could move around and select their own seats, which was different from the 2013–2014 study. The physical flexibility allowed flexibility in building social and academic relationships.

Two types of information were gathered. First students answered questions on their perceptions: (1a) impact of the club and (1b) who is good at math. Second students answered questions on their relationships (2a) who their friends are and (2b) with whom they work. Survey information from 1a and 1b was tabulated. However, the information from 2a and 2b resulted in two directed networks (or graphs) with each node a student and each edge a relationship. Social network analysis was applied to both networks evaluating the centrality of each node (e.g., importance of the student to the network).

Club activities were designed and conducted by Lawrence Smolinsky and his graduate assistant Andrew Alaniz. Observations and surveys were designed and conducted by Eugene Kennedy and his graduate assistant, Liuli Huang.

**Participants**

A motivation of creating this afterschool math club was to engage students who are not traditionally given enrichment opportunities in math, and so there was no requirement that participants have a record of success or previous high levels of interests in mathematics. Students were recruited from one Title I school. Club attendance was voluntary, and recruitment of students involved addressing the priorities of the teachers and administrators, the concerns of the parents, and the reservations of the students (Kennedy & Smolinsky, 2016).

The participants in the present study were 19 middle school students. A total of 21 students participated in the club during the school year, with 19 completing the year. The school setting was a STEM focused charter school that promoted participation in a variety of STEM clubs and activities. Of these 19 participants, there were 6 African American males, 3 African American females, 2 Asian females, 2 Hispanic males, 2 White males, and 1 White female. Two males and 1 female identified as *other* or multi-racial. We conducted statistical comparisons by race or gender using the Mann-Whitney U, Kruskal–Wallis, and Dunn's tests.

**Activities**

While Kennedy and Smolinsky (2016) discussed the intellectual concept and content of the math circle activities, the activities included in the present study were organized with more opportunities for group activities where students could work on their own or with friends and were free to move about. Students' work activities fell into a combination of the following:

1. Work as individuals or in groups. Students were free to move about, ask questions of each other, change partners or groups, or work on their own depending on their preferences.
2. Guided group activity. The leader (mathematician) stood in the front of the class, asked questions, and waited for response. A response might have been followed by the leader requesting the student to come to the board and explain. Students may talk with adjacent students but not walk around the room.
3. Presentations by the leader. The leader would give a presentation that would involve question/response but limited individual work on the part of students.
4. Game playing. Students would play Atari-go in groups of several students.

There was always a juice and snack break during each session. Students were allowed to continue to work and talk with their group or play Atari-go during the break. Students were allowed a voice in the choice of topics and activities. For math circle activities, students chose from alternatives for the coming week or two. If a math circle activity was too difficult or received a poor response, the leader would conclude it at snack time and proceed to Atari-go starting at snack time.

For example, the topic of *Rational and Irrational Numbers as Decimals* took two sessions. Students worked as individuals or groups through examples to show that fractions are represented as repeating decimals. Then, they constructed an argument in a guided group activity to explain repeating decimals. In a second guided group activity, students concluded that repeating decimals are rational numbers by completing an "if and only if" statement. Next, students worked as individuals or in groups to find examples of decimal numbers they could prove are irrational (i.e., patterns of digits that can be fully described but are not repeating decimals, in contrast to $\sqrt{2}$, $e$, or $\pi$ whose digits are not yet fully known).

We allowed students to continue to work when they seemed to enjoy the activity. One activity we planned was for students to examine multiples of 9 and 11 as they related to sums of digits. The leader expected that students might have been familiar with the multiples of 9 from school. That was not the case, and students worked most of the session as individuals or in groups on multiples of 9 and the sum of digits. The leader then continued with a group presentation and guided activity to construct a proof that a number is a multiple of 9 if and only if its sum of digits is too. Discussion of multiples of 11 was not addressed that year.

Atari-go is a beginning version of the game of Go, the oldest board game in the world. Atari-go was popularized by Yasutoshi Yasuda, a professional Go player who used it as a tool for education and socialization (Yasuda, 2002). In the math club, students would play against each other, but Atari-go was also used in math circle activities. Groups or individuals worked on Go puzzles and would come to the front to report on their progress.

**Survey Measures**

Two types of surveys were given: (1) Opinion questions about the effect of the club and who was good at math and (2) selections for friendship and working relationships.

***Perceptions of the Impact of Participating in the Club (ClubImpact)***

Participants were asked to complete a survey which included 5 Likert items related to their perceptions of the impact of participation in the club. Each item was a positive statement about their relationship to math conditioned on causation of "being in the math club." The questions only measured whether students perceived the club as having an impact. Their attitudes might have changed or not. The items were as follows:

- [Q1] Being in the math club has changed my ideas of what I think math is about.
- [Q2] My confidence in my ability to do math has increased as a result of being in the math club.
- [Q3] As a result of being in the club, I now enjoy challenging math problems.
- [Q4] As a result of being in the math club I have a much more positive view of math.
- [Q5] As a result of being in the math club, I want to take more math courses.

The correlations among scores on the items were all positive 0.75 or higher. Cronbach's alpha for this composite was 0.95. Responses were scored 1-5, and the ClubImpact value was calculated by multiplying the average score across items by 10.

Club members were asked to identify other members of the club who were "good at math." The total number of times each student was identified by a fellow club member was tallied and constitutes this scale. It has a range from 0 (no designations) to 18 (identified by all other participants). The correlation between MathStar and ClubImpact was 0.142 with a sigma (2-tailed) of 0.563, indicating that students' perception of benefit from the club was unrelated to others' perception of their math ability.

### *Friendship and working relationships*

Student were surveyed to define connections among participates. First, each student was asked to identify his/her friends in the club. Second, each student was asked to identify club members with whom they had worked on math problems. For each question, students were provided a list of all club members and asked to circle the names of their friends and students with whom they worked on math in a methodology similar to other friendship studies (Coie et al., 1995; Véronneau & Dishion, 2010).

For both friend and work relationships we considered both reciprocated and nonreciprocated selections. While some researchers have required mutual selection for friendship (e.g., Parker et al., 1995) others have also found benefit in nonreciprocated selection (Newcomb & Bagwell, 1995).

## Networks and Social Network Analysis

We construct two directed graphs: One based on friendship relations and one based on work relations. The nodes in both networks are the students. The relationships in each graph were determined by simple surveys completed by the students during interviews.

### *Directed Networks*

A cornerstone of social network analysis has been graphs that depict connections among entities—math club participants in the current study. In the current study, two items on the student survey were used to define connections among participating students. First, each student was asked to identify his/her friends in the club. Second, each student was asked to identify club members with whom they had worked on math problems. The results were assembled into two directed labeled graphs (i.e., friendship and working networks).

Each node in each network was labeled by an individual student, and a directed edge represents a selection by the student at the tail. Selection as a friend or work partner might not be mutual. There may be a directed edge from node A to node B but not from node B to node A. In friendship, student A may be desirous of a relationship with student B. In a study of children's social networks and obesity, it was found that the claim of friendship might not be met with friendship or even just denial, but may even be met with antipathy (de la Haye et al., 2017). In the working network, student A might feel he answered questions or tutored student B, but not believe they worked together. In many studies of middle school collaboration (e.g. Nathan et al., 2007; Piccolo et al., 2008; Webb et al., 2006) tutoring is viewed as collaboration and students may have been assigned to a collaborative group. In working on challenging problems, students may feel otherwise.

***Centrality Measures in Directed Networks***

Social network analysis usually examines which individuals are influential and may impact many others (Sweet, 2016). These measures of influence are called centrality, and there are many different notions of centrality used for various purposes. Centrality is a number assigned to each node and may reflect:

- Degree centrality–the total number of edges attached to the node. In the case of a directed network, there is an "indegree" centrality and "outdegree" centrality.
- Eigenvector centrality–a measure of influence obtained by again totaling the number of edges with the sum then weighted by the importance of influence of the node attached. For directed graphs, there are both "ineigenvector" and "outeigenvector" centralities.

***Undirected Networks***

To form a friendship or work-on-math community, we required that the members had mutual connections. These relationships would be described as undirected graphs or sociograms of a mutual-friendship network and a mutual-working network. An edge was included between student A and student B if and only if A selected B and B selected A. A major notion explored in social network analysis is how individuals are connected to one another through the network. Individuals might belong to sub-groups or identifiable communities within a larger group. These groupings might form based on demographic characteristics or other attributes and might exert a unique and powerful influence on the behavior and views of individuals.

Cliques are an example of a tight-knight group that can exert significant influence on "the fabric of their relationships with others [and] their levels and types of activity" (Adler & Adler, 1995, p. 145). They are identified as the largest subgroups of students in which each student selected every other student (Luce & Perry, 1949). A community structure exists if the members of the club can be grouped into dense connections within a community (internally) and limited connections outside of the community (externally) (Newman, 2004). There are several methods to obtain communities from an undirected graph. Wolfram's Mathematica 11.0 has five algorithms for computing communities for undirected graphs with *modularity* (Newman 2006). The overall network (rather than individual nodes) can also be examined via its density and "assortativity." Density is the proportion (between 0 and 1) of possible connections that are realized. More densely connected networks are associated with productivity and co-construction of knowledge (Cadima et al., 2012) and highly dense networks may promote a normative climate of conformity and limit individuality and creative expressions by members (Maroulis & Gomez, 2008). One can also ask if students who are more popular tend to associate together and if students who are more in demand as mathematical work partners tend to associate together. This property is called assortativity. The assortativity coefficient is a Pearson correlation coefficient and therefore takes values between -1 and 1 (Newman, 2002). Both the mutual-friendship network and mutual-working network were examined for community structure, cliques, density, and assortativity.

## Results

The results from directed and mutual relationships are considered separately.

### Directed Graphs and Relationships Among Club Position and Outcomes

Figures 1and 2 show the directed graphs constructed from the friendship network and the working network. The number labels on the nodes are student IDs.

**Insert Figure 1**

**Insert Figure 2**

The density of the friendship network was higher than the working network. Students selected (and were selected) 7.21 times out of 18 possible as friends and 3.42 time out of 18 as coworkers. Only 39% of friendships are working relationships and 83% of working relationships are friendships. We tested the eigenvalues and impact on gender and race. The working ineigenvalue was higher for girls than boys although the difference was statistically nonsignificant ($p = 0.108$). Among the racial groups, only the friendship ineigenvector had a statistically significant Kruskal–Wallis test. The multiracial and White groups had higher though statistically nonsignificant results after Bonferroni correction using Dunn's test.

**Insert Table 1**

Table 1 presents correlations (Spearman's Rho) among the participant's social position in the club and impressions from the experience (ClubImpact) as well as perceptions their peers hold of them as being good at math (MathStar). Outcomes from the club experience appear to not yield statistically significant correlations with position in the network, except for indegree, which has a positive relationship with peer perceptions of being good at math (MathStar). This implies that students who are perceived as being good at math are more likely to be identified as friends or as popular with peers than others. This is consistent with the resource view of club connections.

Although not statistically significant, students frequently identified as friends by others tended to report lower levels of impact from the club experience than others as evidenced by the negative correlation between indegree and club impact. Perhaps their popularity made the club experience less important as a social experience.

**Insert Table 2**

Table 2 presents correlations between club position and outcome measures when network position is defined by having club participants identify others with whom they have worked. Several of the centrality measures had statistically significant relationships with ratings by peers as being good at math. Additionally, the outeigenvector had a positive and statistically significant relationship with overall assessment of the club experiences. Reporting working with prominent students was correlated with ClubImpact. It is interesting that prominence itself (i.e., being sought as a work partner) was uncorrelated with ClubImpact (Table 2). Seeking work relationships (i.e., viewing oneself in math relationships) was related to gains in attitude while being sought after was not.

Being perceived by fellow students as good at math and a problem solver was positively correlated with selection as a work partner (indegree). This positive correlation continued to hold when prominence in the club was included in the selection (ineigenvector).

**Insert Table 3**

Table 3 presents correlations among selected centrality measures for the two directed graphs. As one might expect, students receiving designations as friends from their peers and students being selected as working partners in mathematics are positively correlated. In either network, the in-centralities represent being selected by others or being desired for an activity. There were positive though statistically nonsignificant correlations between the in-centrality measures in the two networks. These observations reinforce the picture of the club as one in which social status was associated with perceptions of being good at math. Individuals perceived as being good at math appear to be the most valued friends and to hold a place of prominence.

Less expected are the negative correlations between the in-centralities and out-centralities. These eight correlations were all negative with two being statistical significant. If one views the out-centrality measures in terms of only favorable adjectives like gregarious, confident, and outgoing, this result seems paradoxical. However, the out-centrality measures may also reflect students being needy, solicitous, or seeking friends or help in math. Suppose a student has a large outdegree and small indegree in the working network. He or she claimed more working relationships than the partners he/she selected. It may be because the named student viewed the relationship as a source of getting questions answered or receiving tutoring and not a mutual working relationship. The student may also be perceived in negative terms for reasons unrelated to math. Nevertheless, the out-centrality measures were positively correlated with the club experience as seen in Table 2.

***Undirected Graphs and Communities***

Undirected graphs can be constructed from the mutual relationships between students (see Figures 3 and 4). The density of Fig. 3 shows that 23.4% of all possible mutual friendships were actually realized. The density of mutual working relationships (see Figure 4) was 0.076. The lower density of working relationships was expected as we believed that working relationships would be formed from friendships or turn into friendships. Student relationships can be complex. For example, Students 7, 8, and 9 formed a working triad, but have no pairwise mutual friendships, and students 6, 4, are 20 were a friendship triad but only a working chain. A mutual working relationship does make a mutual friendship relationship likely but not a certainty. If a pair of students worked together, then there was a 69% likelihood they were friends. We observe that there are four 4-student friendship cliques and five 3-student friendship cliques but none of them translated into coherent working relationships.

**Insert Figure 3**

**Insert Figure 4**

The assortativity is negative for both mutual friendship -0.205 and mutual working -0.204. This is atypical for social networks where a positive correlation is common (Newman 2002). For example, coauthorship networks in physics, biology, and mathematics were positive for assortativity (Newman 2002, see Table 1). The negative correlation means it is more likely for popular students to mix with less popular students and students who work with a several students to work with students who do not. It seems a healthy development in a math club and consistent with the resource view of a math club. Dijkstra et al.'s (2013) study of middle school friendship supported that students would try to befriend higher status adolescents and distance themselves from lower status ones. Negative assortativity does not support that phenomenon in the mutual-friendship network, but the lower status of students' higher outeigenvector centrality may point partially in that direction.

Note the 5 isolates in Figure 4, namely students 2, 15, 16, 18, and 19. However, four of these students (2, 15, 16, and 19) seem well connected in friendships in Figure 3. Similarly, Student 1 had only one friend in Figure 1., but that friend connected him or her to a working group in Figure 2. Each of the isolates 2, 15, and 16 are part of some 4-student friendship cliques. Working in isolation does not mean one is socially isolated.

The mutual-friendship graph (see Figure 3) was examined for communities using the modularity algorithm. The classification into communities (FriendGrps) is shown in Figure 3. Two of the FriendGrps formed from the mutual-friendship graph have a good mix of gender and race based compared to the distribution of the whole club (see Table 4). FriendGrp 3 is an exception. It consisted of 5 African-American students and 1 Hispanic student. Each group is one (or a fraction of one) female short of gender equality.

**Insert Table 4**

The mutual-working graph was examined for communities (WrkGrps) using modularity and the result was the connected components (see Figure 4). Call the first three rows WrkGrps 1, 2, and 3. For discussion, the five isolates are labeled group 4. The demographics are presented in Table 5. The composition of the 4 groups are similar with WrkGrp 2 and 3 balanced. The gender balance of the working and friendship groups seems to be an exception to "the near universality of gender segregation in middle childhood and early adolescence" (Gest et al., 2007, p. 43).

**Insert Table 5**

One interrelationship is particularly interesting. FriendGrp 1 consists of one female and 4 of the 5 working isolates who worked independently on their own. It may be that the group formed from students who were mathematically strong and worked independently. However, on another survey, only one of these 5 students answered that he/she agreed or strongly agreed with the statement "I like to work on math problems that make me think" compared to 37% for the club. Table 8 shows that this student ranked the lowest of FriendGrps in MathStar. It seems more likely a social phenomenon or work habits brought the group together rather than mathematical interest or ability.

**Insert Table 6**

**Insert Table 7**

In examining the purely graph theoretic measures for WrkGrps (see Table 6), an unexpected item presented itself. Betweenness is a measure on connected graphs (here a connected directed graph) of the nodes through which the number of direct lines of communication through a node may be counted. We had not expected betweenness to have had relevance for a small club in which direct communication can easily occur. However, WrkGrp 2 had a large average betweenness indicating that many members act as a bridge to other club members. This group, on average, had notably larger outward measures that may contribute to betweenness. However, this group did not have the largest indegree average, but did have the largest ineigenvector average, which includes the prestige, or importance, factor. Examining WrkGrp 2 in the composite survey measures (see Table 7) indicated that this group had the highest average, which indicates level of respect for the group's mathematical ability by the club and a high positive attitude toward the club activity. This potentially gives WrkGrp 2 a central role in the club, but does it translate in social importance? WrkGrp 2 has a high average outeigenvector but not a high average ineigenvector. WrkGrp 2 does not seem to be made up of particularly high social status individuals, but rather the group members are part of a high-status group discussed below (FriendGrp 2). We also note this group is more balanced than others with respect to gender (see Table 5).

**Insert Table 8**

**Insert Table 9**

The five isolates previously discussed were not being perceived as strong in math and largely forming a single friendship group (FriendGrp 1). However, these isolates did not seem socially marginalized. Both the friendship average ineigenvector centrality and average outeigenvector centrality for the isolates are near the grand mean (see Table 6). This trend was also true for indegree and outdegree (not shown). The isolates ranked in the middle (2 and 3) for club impact. While they worked alone, they do not seem socially ostracized nor have bad attitudes.

Turning to the friendship groups, FriendGrp 2 presents itself with the highest composite measures (see Table 8) and, in fact, 4 of the 5 members of WrkGrp 2 are in FriendGrp 2. The absent member of WrkGrp 2 (Student 7) was not an ostracized student but a very popular student with the second highest ineigenvector centrality of the club. FriendGrp 2 had the largest (friendship) betweenness. Hence, WrkGrp 2 was high in both working and friendship betweenness and has central role both mathematically and socially. However, FriendGrp 2 was not the most sought after or named as friends (see Table 9). That top social distinction went to FriendGrp 3, which also had a greater mathematical ineigenvector average. So members of FriendGrp 3 were more sought after as friends and coworkers. FriendGrp 2 had the largest Working outeigenvector, indicating that were being helpful or solicitous whereas their friendship outeigenvector was second of three, indicating they were not the most outgoing or solicitous of friendship. Being desired as a friend and being thought mathematically talented were closely related in this math club.

## Discussion

Previous research suggests that afterschool academic clubs can have a positive impact on student attitudes and achievements; however, there has been relatively less focus in the extant research on the within-club variability of these outcomes for participants. The present study examined the importance of friendship and collaboration in an afterschool math club. Results show that the club was highly connected with regard to participants' friendships, perhaps reflecting the experience of a voluntary extracurricular activity.

### Research Questions

*To what extent and in what ways are students' social relationships dependent upon their working relationships?* Results indicated an interesting cross relationship between the two networks. Naming more working partners, which may mean seeking or being solicitous of help, was negatively correlated with being popular. Naming more friends, which may mean seeking or being solicitous of friends, was negatively correlated with forming working relationships. However, these negative relationships between seeking behavior in one network and popularity in the other network are not due to popular students isolating themselves within each network. This surprising fact was indicated by the negative assortativity, which is not typical within social network phenomena.

Mathematical talent as perceived by participants was rewarded in the social structure of the club. Popularity was positively correlated with being thought mathematically talented and with being named or sought as a working partner. In return, those students perceived as mathematically talented were outgoing and helpful. Seeking or naming more working partners is positively correlated with being thought mathematically talented.

*To what extent is a given student's position in the social structure of the club related to club experiences and/or characteristics of the club members?* It was striking that believing one was in working relationships was strongly tied to believing the club had a positive impact but had a strong negative correlation with being named a friend. Girls appeared more valued as work partners, but the difference between boys and girls was not statistically significant. Likewise, multiracial and White students appeared to be the more valued as friends, but the difference between these groups and other ethnic groups was not statistically significant following Bonferroni correction.

Individuals who worked with mathematically prominent members of the club tended to report higher levels of impact from the experience. This tendency was true for individuals who worked with socially prominent members too, but to a lesser extent. Although the working isolate students did not work on problems with others, they were not social isolates but were a below average group in terms of attitude. It is also possible they were the most talented students in the club but are unrecognized by their peers as such because of their working isolation.

*Are there distinct subgroups or communities within the club? If so, what are their characteristics?* The emergence of both friendship and working communities was delineated in the results. We noted the emergence of a working mathematics community (WrkGrp 2) central to the working and social function of the club. This community scored the highest on club members' opinion of mathematical talent and positive value of the club. They were not the most desired as friends, but were included in the largest friendship group (FriendGrp 2) with more popular students considering them friends. WrkGrp 2 also acted as the go-betweens for mathematical problem solving and knowledge. Hence this core of students based on the mathematical club experience expanded into a larger social group and became a peer resource for the club. The friendship group (FriendGrp 3) that reported the least impact was also the group that named the fewest friends (i.e., they were considered the least outgoing). They were frequently named as friends and appeared to be aloof but popular.

**Implications and Limitations**

One point of methodology we believe is significant is the comparison of two social networks. The examination of a single network is common, but there may be multiple phenomena each described by a distinct network as is done in this study. The methods (comparing two networks) and results of the present study point to the potential of social network analysis as a tool to obtain systematic results on student clubs, out-of-school programs, and small classrooms. The findings also point to the need to critically examine social roles, normative climates and status hierarchies in these settings.

**Limitations**

The greatest limitation of the study is its size. The club size was 19 students, which is typical of extra-curricular math club but small from a statistical standpoint. It points to some intriguing results but needs further study. A second limitation is the firm establishment of the reliability and validity of the ClubImpact questionnaire. The method of constructing the networks with student nominations is well established. Also, the specific questions for ClubImpact are direct and pass face validity but should be studied on a larger sample for criterion-related validity (DeVellis, 2012). Furthermore, the reliability check was on a small and dependent sample and should be expanded for a larger study.

## Conclusion

The findings point to the need to critically examine social roles, normative climates, and status hierarchies. The notion of separating friendships and working or academic relationships may be extended to classroom studies in middle schools and high schools.

The immediate extension of this research for the middle school math clubs, is to increase the size of the sample. To increase sample size and thereby the statistical power, one could increase the number of math clubs rather than increase the size the club which might be hard to do. Well-defined network theoretic algorithms, measures, and procedures can allow the results from distinct math clubs to be compared and to draw firm conclusions. These procedures could also be applied in the study of networks of working and social relationships in classrooms.

## Acknowledgements

This work was partially supported by the National Science Foundation (VIGRE award #0739382) and the Louisiana Department of Education 21st Century Community Learning

Centers (award #678PUR-21). The authors are grateful to the two anonymous reviews whose suggestions improved the organization of the article.

# Tables

## Table 1

*Correlations Among Club Position and Outcomes: Friendship Designations*

| | | ClubImpact | MathStar |
|---|---|---|---|
| Spearman's rho | Indegree | -.212 | .443 |
| | Sig. (2-tailed) | .384 | .057 |
| | Ineigenvector | -.249 | .196 |
| | Sig. (2-tailed) | .305 | .421 |
| | Betweenness | .318 | .314 |
| | Sig. (2-tailed) | .184 | .190 |

## Table 2

*Working Network Correlations*

| | | ClubImpact | MathStar |
|---|---|---|---|
| Spearman's rho | Indegree | -.104 | .551 |
| | Sig. (2-tailed) | .672 | .014 |
| | Outdegree | .377 | .338 |
| | Sig. (2-tailed) | .112 | .157 |
| | Outeigenvector | .544 | .258 |
| | Sig. (2-tailed) | .016 | .287 |
| | Ineigenvector | -.023 | .566 |
| | Sig. (2-tailed) | .926 | .012 |
| | Betweenness | .223 | .479 |
| | Sig. (2-tailed) | .360 | .038 |

## Table 3

*Correlation Among Club Position: Working and Friendship Designations*

| Friendship / Working | Indegree | Outdegree | Outeigenvector | Ineigenvector |
|---|---|---|---|---|
| Indegree | .452 | -.218 | -.253 | .336 |
| Sig. (2-tailed) | .052 | .369 | .295 | .160 |
| Outdegree | -.249 | .385 | .330 | -.366 |
| Sig. (2-tailed) | .304 | .103 | .167 | .124 |
| Outeigenvector | -.466 | .278 | .207 | -.589 |
| Sig. (2-tailed) | .044 | .250 | .395 | .008 |
| Ineigenvector | .444 | -.085 | -.114 | .331 |
| Sig. (2-tailed) | .057 | .731 | .641 | .166 |

**Table 4**

*Subgroup Communities: Demographics for Mutual Friendship*

| FriendGrp | Count | STD_ID | Proportion Female | Race Groups |
|---|---|---|---|---|
| 1 | 5 | 2,5,15,18,19 | .40 | 4 |
| 2 | 8 | 1,8,9,10,11,12,13,21 | .375 | 4 |
| 3 | 6 | 3,4,6,7,16,20 | .333 | 2 |

**Table 5**

*Subgroup Communities: Demographics for Mutual Working*

| WrkGrp | Count | ID | Proportion Female | Race Groups |
|---|---|---|---|---|
| 1 | 4 | 1,11,12,21 | .25 | 3 |
| 2 | 5 | 7,8,9,10,13 | .40 | 3 |
| 3 | 5 | 3,4,5,6,20 | .60 | 3 |
| Isolates | 5 | 2,15,16,18,19 | .20 | 3 |

**Table 6**

*Subgroup Communities: Mutual Working*

| WrkGrp | | Friendship average Outeigenvector | Friendship average Ineigenvector | Working average Indegree | Working average Outdegree | Working average Ineigenvector | Working average Outeigenvector | Working average Betweenness |
|---|---|---|---|---|---|---|---|---|
| 1 | Mean | 0.0402 | 0.0252 | 3.25 | 3 | 0.049 | 0.063 | 15.071 |
| N=4 | Std. Dev | 0.0216 | 0.0091 | 0.5 | 2.708 | 0.005 | 0.037 | 17.024 |
| 2 | Mean | 0.0700 | 0.0449 | 4 | 6.4 | 0.068 | 0.113 | 50.021 |
| N=5 | Std. Dev | 0.0380 | 0.0256 | 1.871 | 3.286 | 0.025 | 0.037 | 33.271 |
| 3 | Mean | 0.0447 | 0.0776 | 4.4 | 2.8 | 0.059 | 0.019 | 26.293 |
| N=5 | Std. Dev | 0.0184 | 0.0135 | 2.302 | 1.304 | 0.022 | 0.009 | 38.876 |
| Isolates | Mean | 0.0531 | 0.0573 | 2 | 1.4 | 0.034 | 0.018 | 0.229 |
| N=5 | Std. Dev | 0.0104 | 0.0302 | 1.225 | 2.074 | 0.020 | 0.026 | 0.511 |

**Table 7**

*Subgroup Communities: Mutual Working*

| WrkGrp | | ClubImpact | MathStar |
|---|---|---|---|
| 1 | Mean | 36 | 2.75 |
| | Std. Dev | 7.303 | 1.5 |
| 2 | Mean | 43.6 | 3.6 |
| | Std. Dev | 7.266 | 0.894 |
| 3 | Mean | 28 | 3 |
| | Std. Dev | 10.1 | 1.414 |
| Isolates | Mean | 31 | 1.4 |
| | Std. Dev | 5.568 | 0.548 |

**Table 8**

*Subgroup Communities: Mutual-friendship*

| FriendGrp | | ClubImpact | MathStar |
|---|---|---|---|
| 1 | Mean | 34.2 | 1.8 |
| | Std. Dev | 5.119 | 1.304 |
| 2 | Mean | 39.5 | 3.13 |
| | Std. Dev | 8.124 | 1.246 |
| 3 | Mean | 28.33 | 2.83 |
| | Std. Dev | 10.985 | 1.329 |

**Table 9**

*Subgroup Communities: Demographics and Social Structure based for Mutual-friendship*

| FriendGrp | | Working Outeigenvector | Working Ineigenvector | Friendship Indegree | Friendship Outdegree | Friendship Outeigenvector | Friendship Ineigenvector | Friendship Betweenness |
|---|---|---|---|---|---|---|---|---|
| 1 | Mean | 0.021 | 0.023 | 7.4 | 8 | 0.058 | 0.059 | 11.771 |
| N=5 | Std. Dev | 0.022 | 0.016 | 4.45 | 1.732 | 0.014 | 0.034 | 9.231 |
| 2 | Mean | 0.090 | 0.061 | 5.25 | 7.38 | 0.053 | 0.03 | 15.636 |
| N=8 | Std. Dev | 0.044 | 0.018 | 1.982 | 5.755 | 0.038 | 0.014 | 18.152 |
| 3 | Mean | 0.030 | 0.066 | 9.67 | 6.33 | 0.047 | 0.078 | 11.677 |
| N=6 | Std. Dev | 0.032 | 0.038 | 2.16 | 3.266 | 0.023 | 0.015 | 15.591 |

# Figures

## Figure 1

*Social Network Graph Based on Friendship Designations*

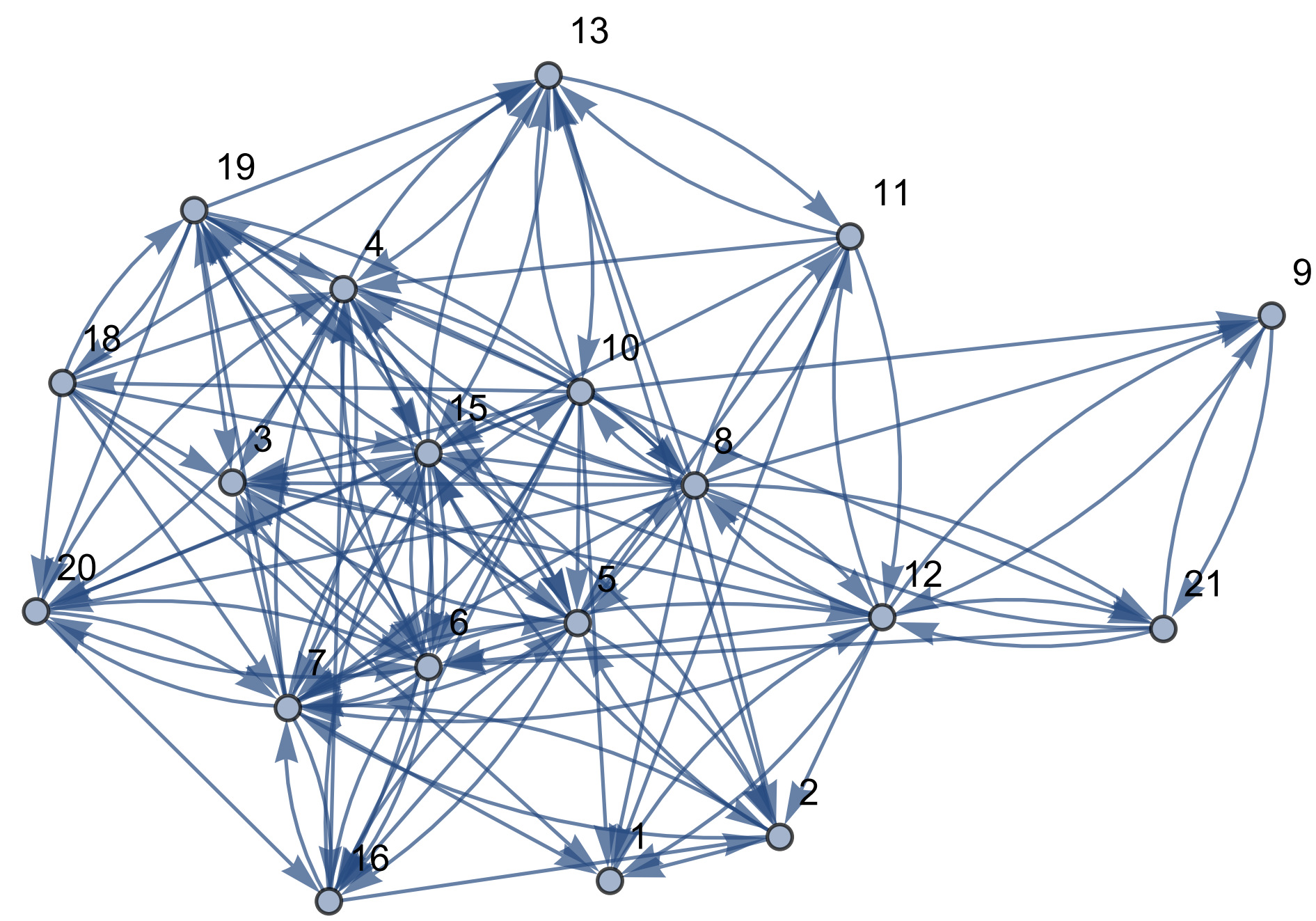

**Figure 2**

*Social Network Graph Based on Working Together on Math Problems Designations*

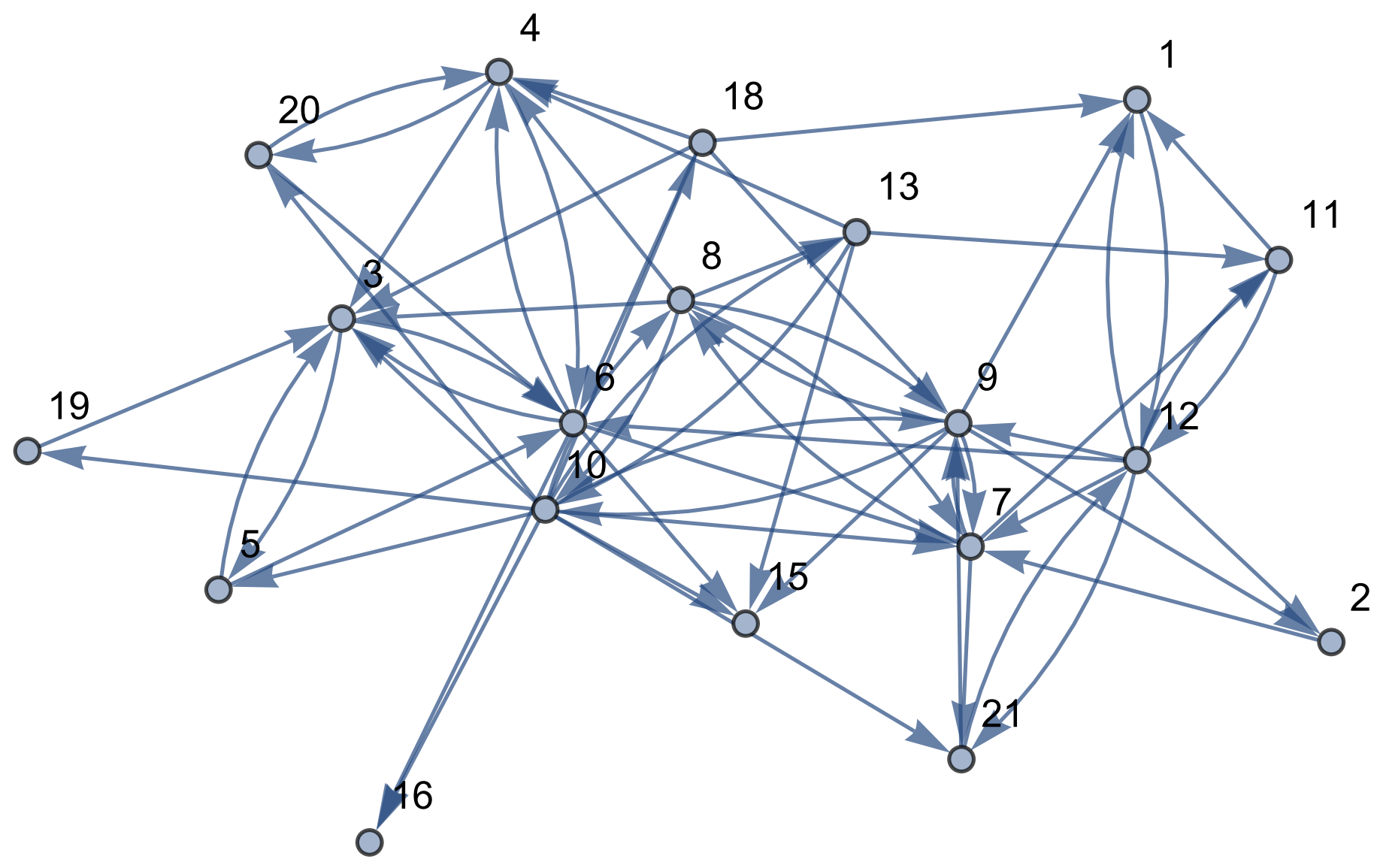

**Figure 3**

*Social Network Undirected Graph Based on Friendship Designations*

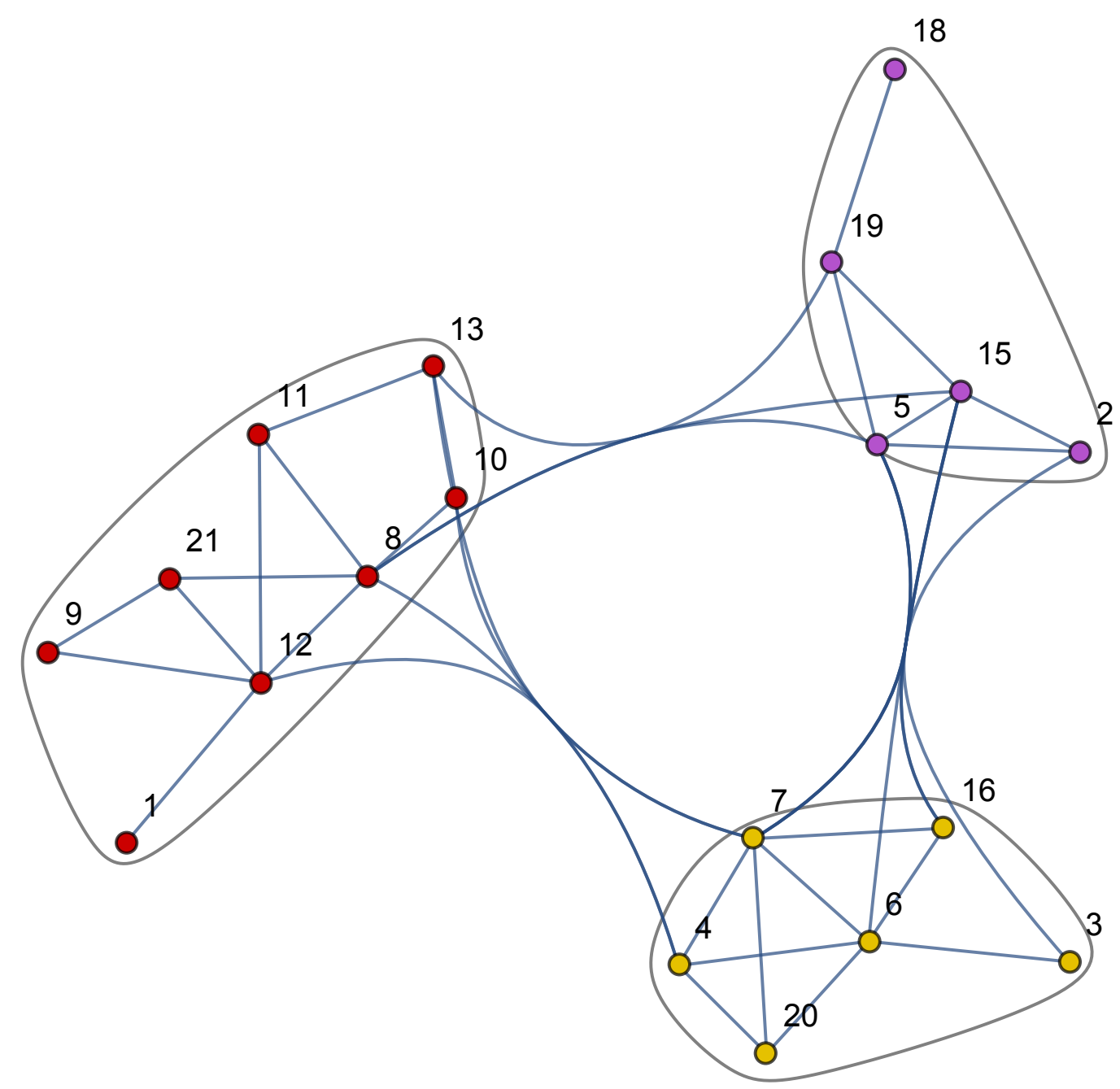

*Note*. The three communities are based on modality.

**Figure 4**

*Social Network undirected graph based on working designations*

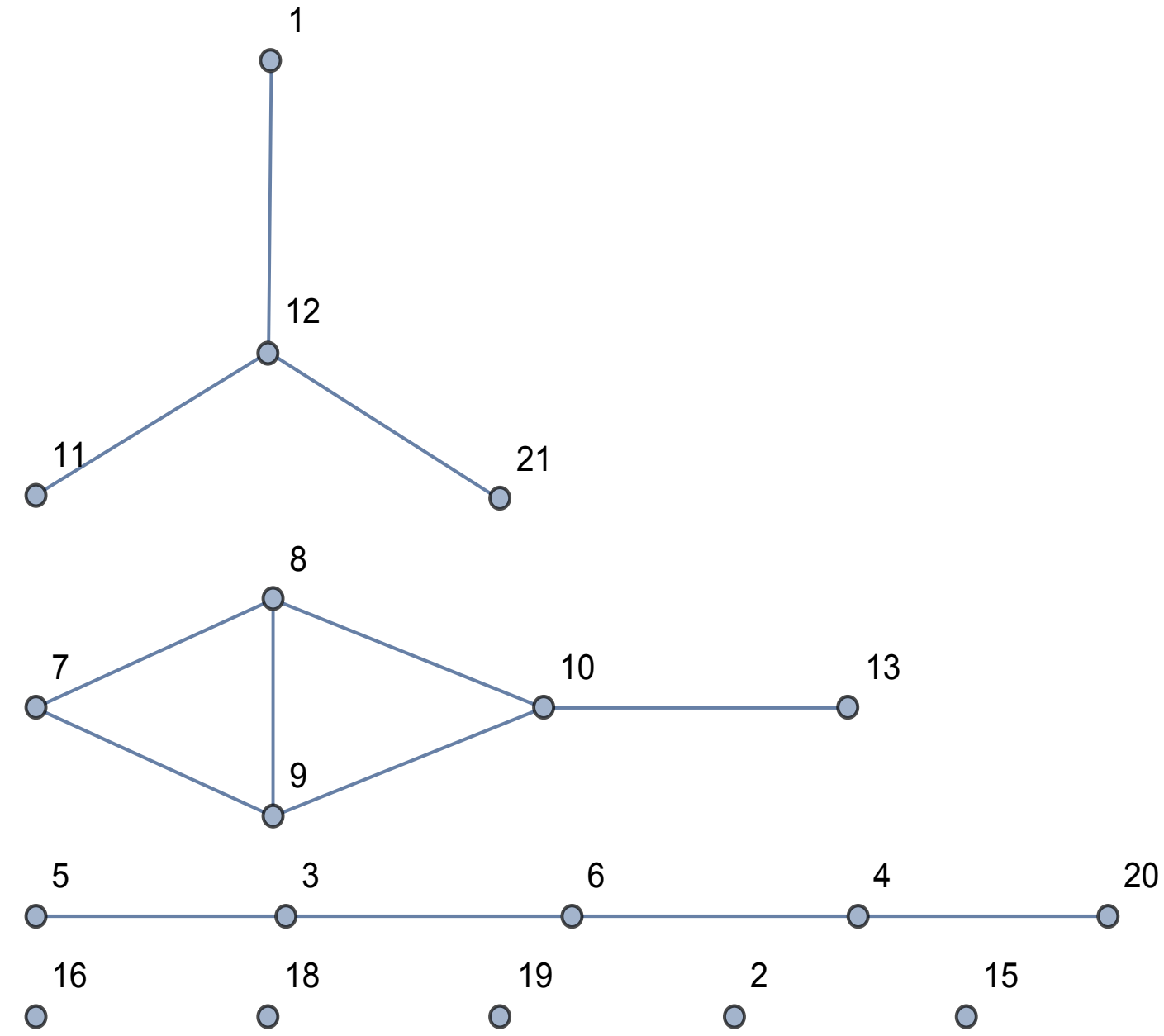

*Note*. The three multimode components are communities based on modality